\documentclass[a4paper,11pt,english]{smfart}
\usepackage{amsfonts}
\usepackage{amsmath,amsthm}
\usepackage{latexsym}
\usepackage{array}
\usepackage{amssymb}
\usepackage{media9}
\usepackage{graphicx}
\usepackage{enumerate}
\usepackage[english]{babel}

\usepackage{color}
\definecolor{marin}{rgb}   {0.,   0.3,   0.7} 
\definecolor{rouge}{rgb}   {0.8,   0.,   0.} 
\definecolor{sepia}{rgb}   {0.8,   0.5,   0.} 
\usepackage[colorlinks,citecolor=marin,linkcolor=rouge,
            bookmarksopen,
            bookmarksnumbered
           ]{hyperref}

\newtheorem{lemma}{Lemma}[section]

\newtheorem{remark}[lemma]{Remark}
\newtheorem{example}[lemma]{Example}

\newtheorem{notation}[lemma]{Notation}
\newtheorem{definition}[lemma]{Definition}
\newtheorem{conclusion}[lemma]{Conclusion}
\numberwithin{equation}{section}
\newcommand{\QED}{\mbox{}\hfill \raisebox{-0.2pt}{\rule{5.6pt}{6pt}\rule{0pt}{0pt}} 
          \medskip\par}

\newcommand{\e}{\mathrm{e}}

\newcommand{\dd}{\,\mathrm{d}}

\newcommand{\Hc}{\mathcal{H}}

\newcommand{\R}{\mathbb{R}}

\newcommand{\Gc}{\mathcal{G}}

\newcommand{\Z}{\mathbb{Z}}

\newcommand{\Kc}{\mathcal{K}}

\let\oldmarginpar\marginpar
\renewcommand\marginpar[1]{\-\oldmarginpar[\raggedleft\footnotesize #1]%
{\raggedright\footnotesize #1}}

\author{Nicolas Crouseilles}
\address{INRIA \& IRMAR  \\
 F-35042 Rennes, France. } 
\email{nicolas.crouseilles@inria.fr}

\author{Lukas Einkemmer}
\address{University of Innsbruck  \\
Technikerstra\ss e 19a, 
A-6020 Innsbruck, Austria. } 
\email{lukas.einkemmer@uibk.ac.at}

\author{Erwan Faou}
\address{INRIA \& ENS Cachan Bretagne  \\
Avenue Robert Schumann F-35170 Bruz, France. } 
\email{Erwan.Faou@inria.fr}

\title[A Hamiltonian splitting for the Vlasov--Maxwell system]
{Hamiltonian splitting for the Vlasov--Maxwell equations}

\begin{document}

\begin{abstract}
A new splitting is proposed for solving the Vlasov--Maxwell system. This splitting is based on a decomposition of the Hamiltonian of the Vlasov--Maxwell system and allows for the construction of arbitrary high order methods by composition (independent of the specific deterministic method used for the discretization of the phase space). Moreover, we show that for a spectral  
method in space this scheme satisfies Poisson's equation without explicitly solving it. Finally, we present some examples in the context of the time evolution of an electromagnetic 
plasma instability which emphasizes the excellent behavior of the new splitting compared to methods from the literature. 

\end{abstract}

\thanks{This work is supported by the Fonds zur F\"orderung der Wissenschaften (FWF) -- project id: P25346.
}
\thanks{This work was partly supported by the ERC Starting Grant Project GEOPARDI}

\maketitle
\tableofcontents

\section{Introduction}

The time evolution of an ensemble of charged particles in a plasma or a propagating beam is described by the Vlasov equation. This kinetic model describes the plasma response to electromagnetic fields. The unknown of the Vlasov equation is the distribution function $f(t, x, v)$ of the considered plasma system, which depends on the time $t$, the space $x$, and the velocity $v$. The coupling of the Vlasov equation with Maxwell's equations takes into account the self-consistent fields, {\it i.e.} the electric $E$ and magnetic field $B$ generated by the particles in the plasma. 

The so-obtained Vlasov--Maxwell system is nonlinear and thus analytical solutions are 
not available in general. Therefore, numerical simulations have to be conducted in order to study realistic physical phenomena (such as magnetic field generation). 

Historically, Particle In Cell (PIC) methods have been widely used to numerically solve high dimensional Vlasov--Maxwell problems. The main advantage of this class of methods is their low computational cost \cite{pic, pic3}. In these methods, the trajectories of macro-particles are advanced by using the characteristics curves of the Vlasov equation whereas the electromagnetic fields are computed by gathering the charge and current densities of the simulated particles on a grid in physical space. Note that this grid is only up to three dimensional (as the electric and magnetic fields do not depend on velocity) whereas the phase space can consist of up to six dimensions.

More recently, deterministic methods have been developed \cite{cheng, filbet, crous-filbet, duclous}. In this approach, a phase space grid is used so that finite volumes, finite elements, or finite differences can be employed to approximate the differential operators. These methods usually require more memory when high dimensional problems are considered. However, they do not suffer from numerical heating and statistical noise. In most of the literature, a time splitting approach is employed; in such a scheme the transport in the spatial variable is split from the transport in the velocity variable (\cite{pcp, valis, mangeney, duclous, cpbm}). The electromagnetic fields are then advanced by approximating Maxwell's equations using a suitable time integrator. Finally, let us emphasize that it is well known that for both approaches (PIC and deterministic), the Poisson equation needs to be satisfied to machine precision (even if Poisson's equation is not solved directly in the numerical method). This is the charge conservation problem (\cite{valis, chacon, langdon, villasenor, regine, crous-respaud}). 

In this work, we propose new time splitting schemes for the numerical solution of the Vlasov--Maxwell system using the deterministic approach outlined above. These splitting methods enjoy the following properties: 
\begin{itemize}
    \item  they are {\em symplectic}, in the sense that they preserve the mechanical Hamiltonian structure of the equation; 
	\item they can be generalized to arbitrary high order in time by composition (first order for the Lie splitting, second order for Strang splitting, fourth order for the triple Jump method, $\dots$); 
	\item Poisson's equation is exactly satisfied  at the semi-discrete level as well as at the fully discrete level for standard space discretization methods. 
	\item they can be easily combined with arbitrary high order schemes in phase space (finite volume, semi-Lagrangian, spectral, $\dots$).  
\end{itemize}

These splitting methods are based on a decomposition of the Hamiltonian of the Vlasov--Maxwell equations. 
Splitting methods are widely used in the context of systems of ordinary differential equations governed by Hamilton's equations. In the case of velocity independent potentials for example, the application of the Strang splitting algorithm leads to the well known St\"ormer-Verlet method. 
Recently, these methods have been extensively studied in the literature (see, for example, \cite{hlw} for an overview): as they are symplectic, they exhibit excellent long time behavior and energy preservation properties. 

As far as partial differential equations are concerned, splitting methods are also well developed in the context of linear and semi-linear Hamiltonian equation such as the Schr\"odinger equation (see for instance \cite{Lubich08, F11}). However, in the case of strongly nonlinear  equation, where the highest order derivative depends on the solution itself, the formulation as a Hamiltonian system is significantly less common. Nevertheless, it is well known that the Vlasov--Maxwell equations can be considered as an infinite dimensional Hamiltonian system (see \cite{Morrison80, MaWe82}) where the Hamiltonian structure is non canonical, and depends on the solution itself (Poisson structure).  This is of interest, as contrary to the Vlasov--Poisson equation for which the natural splitting between spatial advection and velocity advection is also a Hamiltonian splitting (see \cite{cfm}), it is not clear to decide from the structure of the Vlasov--Maxwell which choice of splitting scheme yields {\em a priori} good conservation properties. Furthermore, as we will discuss in more detail in section \ref{sec:numerical}, the splitting methods for the Vlasov--Maxwell equations constructed in the literature have exclusively focused on two-term splittings (motivated by the two transport terms in space and in velocity) and therefore have used rather ad-hoc procedures in order to obtain methods of second order. 
This also impedes the construction of high order methods by composition (which is relatively straightforward for the Vlasov--Poisson equation, see \cite{cfm}) as the resulting scheme are no longer symmetric. 

In this paper, we outline the construction and numerical validation of these new splitting methods and apply them to the reduced 1+1/2 dimensional Vlasov--Maxwell system. Several numerical tests have 
been conducted and comparisons with two numerical schemes of the literature (\cite{valis, mangeney}) are performed for Landau damping, a Weibel type instability, and a magnetically induced two-stream instability. 

The rest of the paper is organized as follows: first, the Vlasov--Maxwell model and its Hamiltonian structure are recalled. Then, the Hamiltonian splitting is presented and applied to the reduced 1+1/2 dimensional Vlasov--Maxwell system. Finally, numerical experiments are conducted with the new methods and comparisons with classical methods from the literature are performed. 

\section{Hamiltonian structure of the Vlasov--Maxwell system}

We consider the Vlasov--Maxwell system that is satisfied by the electron distribution function  $f=f(t, x, v)$ and the 
electromagnetic fields $(E,B)=(E(t, x), B(t, x))\in\mathbb{R}^3\times \mathbb{R}^3$. 
Here, the spatial variable is denoted by  $x\in X^3$ ($X^3$ being a three dimensional torus),  
 the velocity variable is denoted by $v\in\mathbb{R}^3$, and the time is denoted by $t\geq 0$. 
Using normalized units, the Vlasov--Maxwell system can be written as 
\begin{equation}
\label{eq:vm}
\begin{array}{l}
\partial_t f + v \cdot \nabla_x f + (E  + v \times B ) \cdot \nabla_v f = 0, \\
\partial_t E = \displaystyle \nabla_x \times B - \int_{\R^3} v f(t,x,v) \dd v, \\
\partial_t B = - \nabla_x \times E. 
\end{array}
\end{equation}
The two constraints on the electromagnetic fields $(E, B)$ are given by
 \begin{equation}
\label{eq:vm2}
\nabla_x\cdot  E = \rho(t,x) := \color{black} \int_{\R^3} f(t, x, v) dv -1, \;\; \nabla_x \cdot B=0. 
\end{equation}
Note that if these constraints are satisfied at the initial time, they are satisfied for all times $t>0$. 
Moreover, the total mass is preserved; that is,  
\[ \int_{X^3}\int_{\R^3}f(t, x, v)\dd x \dd v = \int_{X^3}\int_{\R^3}f(0, x, v)\dd x \dd v=1\]
holds true for all time $t>0$. Note, however, that this is not always the case when numerical approximation are considered. 
For a unique solution the Vlasov--Maxwel system has to be supplemented with initial conditions for the distribution function as well as the field variables; that is, we have to specify 
\[ f(t=0, x, v)=f_0(x, v),\quad E(t=0, x)=E_0(x),\quad B(t=0, x)=B_0(x). \]

The Hamiltonian associated with the Vlasov--Maxwell system is given by (see \cite{Morrison80} and \cite{MaWe82})
\begin{eqnarray}
\Hc &=& \frac12 \int_{X^3} |E|^2 \dd x + \frac12 \int_{X^3} |B|^2 \dd x + \frac12 \int_{X^3 \times \R^3} |v|^2 f \dd v \dd x\nonumber\\
\label{ham-split}
&=& \Hc_E + \Hc_B + \Hc_f. 
\end{eqnarray}
The three terms corresponding to electric energy, magnetic energy, and kinetic energy, respectively.  For a given functional $\Kc(f,E,B)$, we denote by $\delta \Kc/\delta f$, $\delta \Kc/\delta E$ and $\delta \Kc/\delta B$ the Fr\'echet derivatives of $\Kc$ with respect to $f$, $E$ and $B$ respectively\footnote{Note that in the absence of domain boundaries, this definition of functional derivative is not ambiguous.}. 
The Poisson bracket of two functionals $\Kc(f,E,B)$ 
and $\Gc(f,E,B)$, is then defined as \color{black}
\begin{align}
\label{eq:poisson}
[\Kc,\Gc] = & \int_{X^3 \times \R^3} f \{Ê\frac{\delta \Kc}{\delta f }, Ê\frac{\delta \Gc}{\delta f } \}Ê\dd x \dd v \\
 \nonumber & + \int_{X^3}  Ê\frac{\delta \Kc}{\delta E } \cdot (\nabla_x \times Ê\frac{\delta \Gc}{\delta B }) - Ê\frac{\delta \Gc}{\delta E } \cdot (\nabla_x \times Ê\frac{\delta \Kc}{\delta B }) \dd x\\
 \nonumber & + \int_{X^3 \times \R^3} Ê\frac{\delta \Kc}{\delta E } \cdot \nabla_v f  \frac{\delta \Gc}{\delta f} -Ê\frac{\delta \Gc}{\delta E } \cdot \nabla_v f  \frac{\delta \Kc}{\delta f} \dd x \dd v\\
 \nonumber & + \int_{X^3 \times \R^3} \frac{\delta \Kc}{ \delta B} \cdot ( \nabla_v f \times v) \frac{\delta \Gc}{\delta f}
 - \frac{\delta \Gc}{ \delta B} \cdot ( \nabla_v f \times v) \frac{\delta \Kc}{\delta f} \dd x \dd v,  
 \end{align}
where for two functions $h(x,v)$ and $k(x,v)$, 
\[ \{h,k\} = \sum_{i=1}^3 
	\left( \frac{\partial h}{\partial x_i}\frac{\partial k}{\partial v_i}	
     -\frac{\partial h}{\partial v_i}\frac{\partial k}{\partial x_i} \right) \]
	 denotes the standard  microcanonical Poisson bracket. With this notation, the Vlasov--Maxwell system \eqref{eq:vm} is equivalent to the equation	 \begin{equation} \label{eq:eHam}
		 \partial_t \Kc = [\Kc,\mathcal{H}] = [\Kc,\mathcal{H}_E] +[\Kc,\mathcal{H}_B] +[\Kc,\mathcal{H}_f],
	 \end{equation}
for any functional  $\Kc(t) = \Kc(f(t),E(t),B(t))$ evaluated along the solution of \eqref{eq:eHam} (note that here, we assume that the solution $(f,E,B)$ and the functional $\Kc$ are smooth enough to ensure the validity of the equations).   

As is outlined in the next section, this formulation provides the basis for the splitting methods proposed in this paper.

\section{Hamiltonian splitting}
In this section, we propose new splitting methods to compute the solution of the Vlasov--Maxwell system. These splitting are based on exact computations of the three parts of the Hamiltonian $ \Hc_E,  \Hc_B$ and $\Hc_f$ given by \eqref{ham-split} respectively. 

Let us start by detailing the equations associated with the different parts of the Hamiltonian. Note that these evolution equations can be derived by plugging in a representation for $\Hc_E$, $\Hc_B$, and $\Hc_f$ into the right hand side of equation \eqref{eq:eHam}.

\subsection{Equations for $\Hc_E$}
The equations associated with the Hamiltonian $\Hc_E$ are given by 
\begin{equation}
\begin{array}{l}
\partial_t f + E(x) \cdot \nabla_v f = 0, \\
\partial_t E = 0, \\
\partial_t B = - \nabla_x \times E.  
\end{array}
\end{equation}
For a given initial data $(f_0,E_0,B_0)$ at time $t=0$, the solution of this system at time $t$ is given explicitly by 
\begin{equation}
\begin{array}{l}
 f(t,x,v) = f_0(x, v - t E_0(x)),  \\
E(t,x) = E_0(x), \\
B(t,x) = B_0(x) - t \, \nabla_x \times E_0(x). 
 \end{array}
\end{equation}

Formally we can write this solution as 
\[ F(t):=(f, E, B)^T(t)=\exp( \Hc_E t)(f_0, E_0, B_0)^T. \]

Moreover, if $E_0$ and $f_0$ satisfy the relation $\nabla_x\cdot E_0=\int_{\R^3}f_0(x, v)\dd v -1$, 
then it holds that $\nabla_x\cdot E(t,x)=\int_{\R^3}f(t,x, v)\dd v -1$; this can be easily seen by the fact that 
$E(t,x)$ is constant in time, and that the transformation $v \mapsto v - t E_0(x)$ preserves the volume. 

Let us further remark that if $\nabla_x\cdot B_0=0$, then this property holds true for any later time as well. This is easily shown by considering the divergence of the last equation
$$ \nabla_x \cdot B(t, x) = \nabla_x \cdot B_0(x) - t \, \nabla_x\cdot \left(\nabla_x \times E_0(x)\right) = \nabla_x \cdot B_0(x) =0. 
$$

The previous relations can be easily carried over to the fully discrete case by using for example a spectral discretization in $x$ and an interpolation procedure in $v$ to compute $f(t,x,v)$ from $f_0(x,v)$. The volume preservation in this case can be easily ensured as the advection is just a translation with constant coefficients. 

\subsection{Equations for $\Hc_B$}

The equations associated with the Hamiltonian $\Hc_B$ are given by 
\begin{equation}
\begin{array}{l}
\partial_t f + (v \times B(x)) \cdot \nabla_v f = 0, \\
\partial_t E = \nabla_x \times B, \\
\partial_t B = 0.  
\end{array}
\end{equation}
Note that in the right hand side of the first equation the following cross product term appears 
$$
v \times B = \begin{pmatrix}
v_2 B_3  - v_3 B_2 \\
v_3 B_1 - v_1 B_3\\
v_1 B_2 - v_2 B_1
\end{pmatrix}
=
 \begin{pmatrix}
0 &  B_3  &- B_2 \\
-B_3 & 0 & B_1\\
B_2 &- B_1 & 0
\end{pmatrix}
v =: \mathbb{J}_B v.  
$$
Two important properties of this term have to be remarked. First, the $k-$th component of $(v\times B)$ does not depend on $v_k$ (for $k=1, 2, 3$). Second, the $3\times 3$ matrix $\mathbb{J}_B$ is constant in time (which follows immediately from the third evolution equation).

For a given initial data $(f_0, E_0,B_0)$ at time $t=0$, 
the solution of this system at time $t$ is therefore given explicitly by  %
\begin{equation}
\begin{array}{l}
f(t,x,v) = f_0(x, \exp( -\mathbb{J}_B t)v), \\ 
E(t,x) = E_0(x) + t \, \nabla_x \times B_0(x),  \\
B(t,x) = B_0(x). 
\end{array}
\end{equation}
As before, we will denote this formally by 
$$F(t):=(f, E, B)^T(t)=\exp( H_B t)(f_0, E_0, B_0)^T.$$ 

Let us remark that if $\nabla_x\cdot E_0=\int_{\R^3}f_0(x, v)\dd v -1$, then 
considering the divergence of the second equation leads to $\nabla_x\cdot E(t, x)=\nabla_x\cdot E_0(x)$. 

Now as $\mathbb{J}_B$ is skew-symmetric, the transformation $v \mapsto \exp( -\mathbb{J}_B t)v$ preserves the volume, and we thus follow that 
\begin{equation}
\label{eq:volpres}
\int_{\R^3}f_0(x, v)\dd v = \int_{\R^3}f_0(x, \exp( -\mathbb{J}_B t)v)\dd v
\end{equation}
which ensures that the Poisson equation is satisfied after one integration step.

Note that in contrast with the previous case, the relation \eqref{eq:volpres} is in general difficult to ensure by standard 3D (or 2D) interpolation at the fully discrete level. To remedy this difficulty, one possibility is to further approximate $f(t,x,v)$ by directional splitting.  
Indeed, due the structure of the cross product, {\it i.e.} since 
$$
v \times B = \begin{pmatrix}
v_2 B_3  - v_3 B_2 \\
v_3 B_1 - v_1 B_3\\
v_1 B_2 - v_2 B_1
\end{pmatrix},
$$
the evolution equation for $f(t,x,v)$ can be written as
\[	\partial_t f + (v_2 B_3  - v_3 B_2)\partial_{v_1}f + (v_3 B_1 - v_1 B_3)\partial_{v_2}f + (v_1 B_2 - v_2 B_1)\partial_{v_3}f = 0, \]
which can be further split into the following three equations
\begin{align*}
	\partial_t f &= -(v_2 B_3  - v_3 B_2)\partial_{v_1} f \\
	\partial_t f &= -(v_3 B_1 - v_1 B_3)\partial_{v_2}f \\
	\partial_t f &= -(v_1 B_2 - v_2 B_1)\partial_{v_3}f.
\end{align*}
Each of these equations can be solved explicitely by a linear advection with constant coefficients. For example the speed of the advection in the $v_1$ direction does only depend on the perpendicular velocities $v_2$ and $v_3$. 
As proposed in \cite{mangeney}, a Strang splitting (or high order splitting) can be performed here. 

\subsection{Equations for $\Hc_f$}

Finally, the equation associated with the Hamiltonian $\Hc_f$ is given by 
\begin{equation}
\begin{array}{l}
\partial_t f + v \cdot \nabla_x f  = 0, \\
\partial_t E =  - \int_{\R^3} v f(t,x,v) \dd v, \\
\partial_t B = 0. 
\end{array}
\end{equation}
For a given initial data $(f_0, E_0,B_0)$ at time $t=0$, the solution at time $t$ is given explicitly by 
\begin{equation}
\begin{array}{l}
f(t,x,v) = f_0(x - tv,v), \\
B(t,x) = B_0(x), \\
E(t,x) = E_0(x) - \displaystyle\int_0^t \int_{\R^3} v f_0(x - sv,v) \dd v \dd s.   \\
 \end{array}
\end{equation}
We will denote this solution formally by 
$$F(t):=(f, E, B)^T(t)=\exp( H_f t)(f_0, E_0, B_0)^T.$$

Let us remark that the Poisson equation is propagated with time. Indeed, 
considering the divergence of the last equation leads to 
\begin{eqnarray*}
\nabla_x \cdot E(t, x)&=& \nabla_x \cdot E_0(x) - \int_0^t \int_{\R^3} v\cdot \nabla_x [ f_0(x - sv,v)] \dd v \dd s \nonumber\\
&=& \int_{\R^3} f_0(x, v)\dd v - 1 + \int_0^t \int_{\R^3} \partial_s [ f_0(x - sv,v)] \dd v \dd s \nonumber\\
&=& \int_{\R^3} f_0(x, v)\dd v - 1 + \int_{\R^3} [ f_0(x - tv,v) - f_0(x, v)] \dd v \nonumber\\
&=&  \int_{\R^3}  f_0(x - tv,v) \dd v - 1\nonumber\\
&=&  \int_{\R^3}  f(t, x,v) \dd v - 1,  
\end{eqnarray*}
which implies that the splitting method proposed here satisfies the charge conservation property. 

The preservation of the charge conservation at the fully discrete level relies on a good approximation of the time integral in the previous formula, in combination with the interpolation procedure in the $v$ variable. 
We will give an example in the next section where $E(t,x)$ is easily computed using Fourier spectral approximations in the $x$ variable. 

\subsection{Hamiltonian splitting}
\label{sec3.4}
The solution $F(t)=(f, E, B)^T(t)$ of the Vlasov--Maxwell system up to first order can be written as (this is the well known Lie splitting)
$$
F(t) = \exp(\Hc_E t)\exp(\Hc_B t)\exp(\Hc_f t) F(0). 
$$
The Strang splitting scheme (which is a second order accurate splitting method) is given by
$$F(t) =  \exp(\Hc_E t/2)\exp(\Hc_B t/2)\exp(\Hc_f t) \exp(H_B t/2)\exp(H_E t/2)F(0).
$$
Higher order methods can be constructed using composition (see  \cite[III.3]{hlw}). For example the well known triple jump scheme (see \cite[II, Example 4.2]{hlw}) based on the Strang splitting yields a method of order four, as the original Strang splitting is symmetric.

\section{Application to the reduced 1 + 1/2 dimensional model}
In this section, we apply the splitting method proposed in this paper to a fully discretized 1+1/2 model of the Vlasov--Maxwell system. This model has been studied in \cite{pcp, cglm} and retains most of the properties of the full Vlasov--Maxwell system. All the numerical simulations in this paper will be conducted in the framework of this reduced model. 

\subsection{Reduced 1+1/2 model}
We consider the phase space $(x_1,v_1,v_2) \in X\times \R^2$, where $X$ is a one-dimensional torus,  and the unknown functions 
$f(t,x_1,v_1,v_2)$, $B(t,x_1)$ and $E(t,x_1) = (E_1,E_2)(t,x_1)$ which are determined by solving the following system of evolution equations

\begin{equation}
\label{eq:vmred}
\begin{array}{l}
\partial_t f + v_1  \partial_{x_1} f + E \cdot \nabla_v f  + B \mathcal{J} v \cdot \nabla_v f = 0,\\
\partial_t B = -  \partial_{x_1} E_2; \\
\partial_t E_2 = -  \partial_{x_1} B - \displaystyle \int_{\R^2} v_2 f(t,x_1,v) \dd v,\\
\partial_t E_1 = - \displaystyle \int_{\R^2} v_1 f(t,x_1,v) \dd v,\\
\end{array}
\end{equation}
where $v=(v_1, v_2)$ and $\mathcal{J}$ denotes the symplectic matrix 
$$
\mathcal{J} = \begin{pmatrix}
0 & 1 \\Ê-1 & 0
\end{pmatrix}.
$$
This reduced system correspond to choosing an initial value of the form 
$$
E(x_1,x_2,x_3) = \begin{pmatrix} E_1(x_1) \\ E_2(x_1) \\ 0 \end{pmatrix}
\quad \mbox{and}\quad B(x_1,x_2,x_3) = \begin{pmatrix} 0 \\ 0 \\ÊB(x_1) \end{pmatrix}
$$
and a $f$ depending on $x_1$ and $(v_1,v_2)$ only. 
Then it can be easily checked that this structure is preserved 
by the flow of \eqref{eq:vm} and that the equations reduce to \eqref{eq:vmred}. 
We refer the reader to \cite{cglm} for more details.

\subsection{Phase space integration}
Consistent with most of the literature, 
we consider the periodic case in the $x_1$ direction. Thus, one natural method in this direction 
is the Fourier spectral approximation. It simplifies the computation of the evolution corresponding 
to $\Hc_f$, and the exact charge conservation can be easily recovered. 
In the evolutions corresponding to $\Hc_E$ and $\Hc_B$,  
we compute an approximation to the transport of $f$ by using a third order finite volume numerical scheme (see \cite{filbet, cms}) whereas the equations for
the electric field $E$ and the magnetic field $B$ are solved in Fourier space. 

Let us consider the evolution corresponding to $\Hc_f$ in more detail. 
By employing the Fourier transform in space (that is in the $x_1$  variable) 
and denoting the corresponding functions by $\hat{f}_k$,  $k \in \Z$, (that is, $\hat{f}_k$ 
denotes the Fourier transform of $f(x)$ and does still depend on $v_1$ and $v_2$) 
gives
$$
\partial_t \hat{f} + v_1 i k \hat{f} = 0, \;\;\;\; \partial_t \hat{E}_1 = -\int_{\R^2} v_1 \hat{f} \dd v. 
$$
The first equation can be solved exactly in time, given the initial condition  $\hat{f}(t^n)$: 
$$
\forall\, k \in \Z, \quad 
\hat{f}(t) = \hat{f}(t^n) \exp(- i v_1 k (t-t^n)), 
$$
which inserted into the second equation (which we have integrated in time from $t^n$ to $t^{n+1}$) yields the approximation  
\begin{eqnarray*}
\hat{E}_1^{n+1} &=& \hat{E}_1^n - \int_{t^n}^{t^{n+1}} \int_{\R^2} v_1 \hat{f}(t^n) \exp(- i v_1 k (t-t^n)) \dd v \dd t \nonumber\\
  &=&  \hat{E}_1^n - \int_{\R^2} v_1 \hat{f}(t^n) \int_{t^n}^{t^{n+1}} \exp(- i v_1 k (t-t^n)) \dd t \dd v  \nonumber\\
    &=&  \hat{E}_1^n - \int_{\R^2} v_1 \hat{f}(t^n) \left[   \frac{-1}{i k v_1} (\exp(-ikv_1 \Delta t)-1 )  \right] \dd v.   \nonumber\\
\end{eqnarray*}
For the long time properties of the numerical scheme it is vital that the total
charge is conserved exactly and that PoissonÕs equation is satisfied. 
The charge is not part of the formulation of our numerical algorithm, but 
it can be shown that our scheme satisfied the charge conservation property. Indeed, assuming that Poisson's equation is satisfied initially, we can prove that 
 it is satisfied for all times. 

Assuming that Poisson equation is satisfied at time $t^n$, that is $ik\hat{E_1}^n =\hat{\rho}^n$, where $\hat{\rho}^n$ is an approximation of the charge density $\rho(t,x)$ in Fourier, see \eqref{eq:vm2}, 
we get from the last equation 
\begin{eqnarray*}
\hat{E}_1^{n+1} &=&  \hat{E}_1^n - \int_{\R^2} v_1 \hat{f}(t^n) \left[   \frac{-1}{i k v_1} (\exp(-ikv_1 \Delta t)-1)   \right] \dd v  \nonumber\\
&=& \frac{\hat{\rho}^n}{ik} + \frac{1}{ik}   \int_{\R^2}  \hat{f}(t^n) (\exp(-ikv_1 \Delta t)-1)  \dd v  \nonumber\\
&=& \frac{1}{ik} \left[ \hat{\rho}^n  +  \int_{\R^2}  \hat{f}(t^n)  \exp(-ikv_1 \Delta t) \dd v - \hat{\rho}^n -1\right] \nonumber\\
&=& \frac{1}{ik} \hat{\rho}^{n+1},  
\end{eqnarray*}
which shows that Poisson's equation is also true at time $t^{n+1}$,  
assuming a sufficient accuracy in the discretization of the integral in $v$. 
Note that $k=0$ does not pose any additional difficulties as the average of the electric field is imposed to be 
zero by total mass conservation. This is a necessary condition to impose as periodic boundary conditions 
do not give a unique solution for the electric field. 

Once the three steps are solved, the global time integration is given in subsection \ref{sec3.4}. 

\section{Numerical experiments}
\label{sec:numerical}
In this section, the Hamiltonian splitting is numerically studied by solving the 1+1/2 dimensional reduced Vlasov--Maxwell system introduced in the previous section. We also compared our numerical scheme to two methods of the literature:  
VALIS (from \cite{valis}, recalled in subsection \ref{algo-valis}) and Mangeney 
(from \cite{mangeney}, recalled in \ref{algo-mangeney}).  

\subsection{Landau damping}
Landau damping is a classic test case for the Vlasov--Poisson equations. 
In this context a number of numerical schemes have been proposed (see e.g. \cite{filbet,cms,cfm}). 
However, even though the equation we try to solve does not explicitly determine the electric field from the charge density (as is usually done for methods constructed in the context of the Vlasov--Poisson equations) it clearly reduces to the Vlasov--Poisson equation if no magnetic effects are present in a given configuration. Compared to the simulation discussed in sections \ref{sec:tsi} and \ref{sec:tt} this problem does not rely on magnetic effects. Thus, it serves as a first, but rather limited, consistency check for the numerical schemes considered in this paper. However, the charge conservation property 
is crucial in this test; if it is not satisfied exactly, the solution of the Vlasov--Poisson system is not recovered.  

For this test, the initial value for the distribution function $f$ is given by
\[
f_0(x_1,v_1,v_2) = \frac{1}{2\pi}\e^{-\frac{1}{2}\vert v \vert^2}\left(1+\alpha \cos kx_1 \right), 
\]
where $x_1\in X=[0, 2\pi/k]$, $v\in \R^2$ and we have chosen $\alpha=0.5$ and $k=0.4$, a configuration usually called strong (or nonlinear) Landau damping. 
The electric field, at the initial time, is determined by Poisson's equation, {\it i.e.}
$$
\partial_{x_1} E_1(t=0, x_1) = \int_{\R^2} f_0(x_1, v) \dd v-1, 
$$
whereas a vanishing magnetic and perpendicular electric field ({\it i.e.}~\mbox{$E_2$}) is prescribed in order to completely specify the initial value.
	
The quantity of interest is the potential energy stored in the electric field
$$
{\mathcal E}_{pot}(t) = \frac{1}{2}\int_X |E(t, x_1)|^2 \dd x_1, 
$$
and the total energy (sum of the potential and kinetic energy) 
$$
{\mathcal E}_{tot}(t) ={\mathcal E}_{pot}(t)+{\mathcal E}_{kin}(t), $$
where
$$ {\mathcal E}_{kin}(t)=\frac{1}{2}\int_X\int_{\R^2}|v|^2f(t, x_1, v) \dd v \dd x_1.$$

Note that the total energy is conserved with time.

\begin{figure}
	\includegraphics[width=12cm]{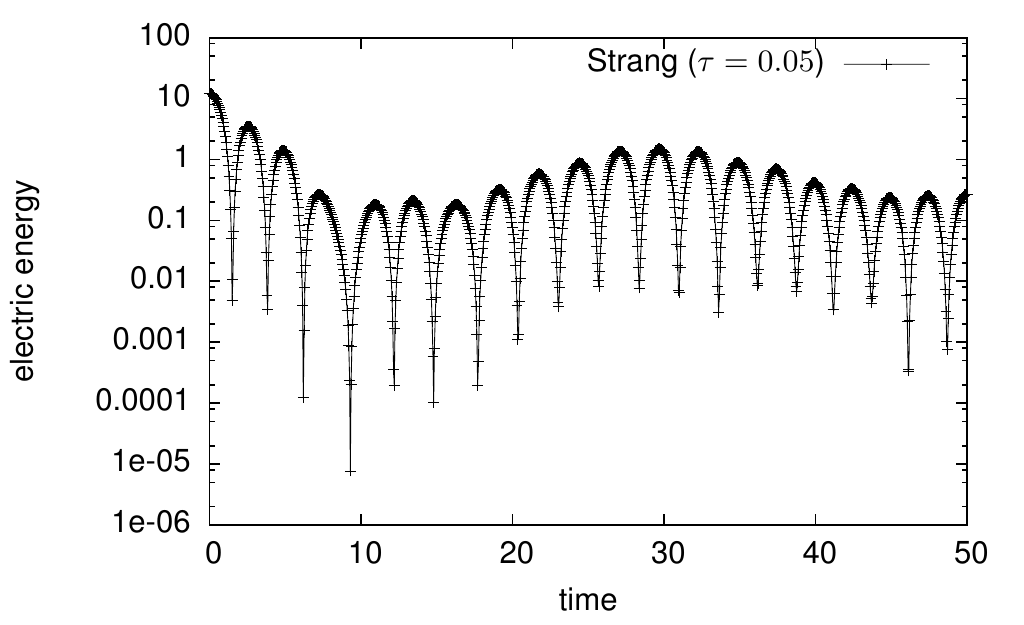}
	\includegraphics[width=12cm]{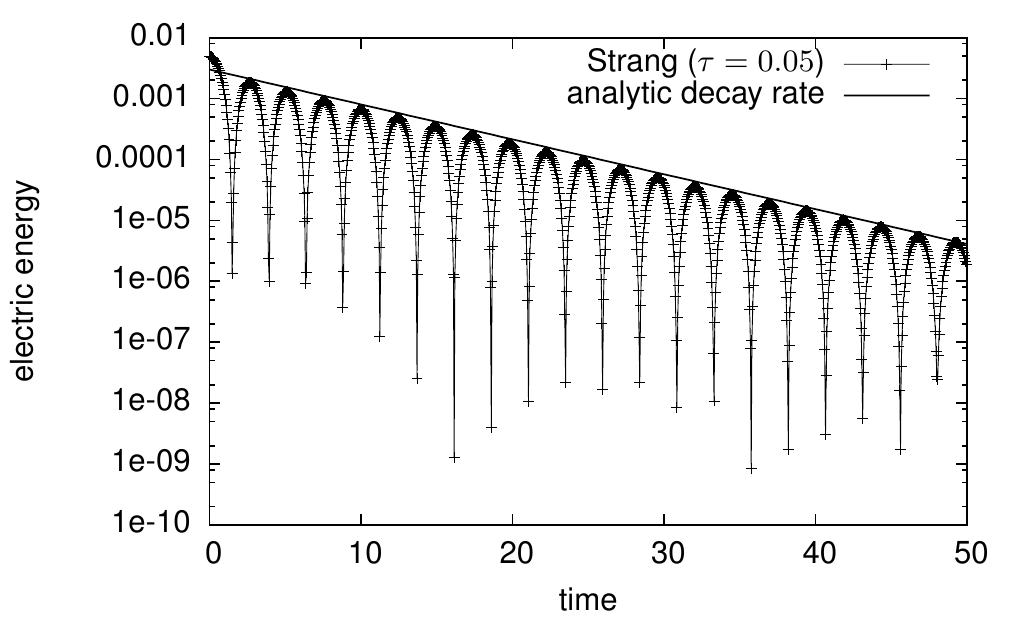}
	\caption{Time evolution of the electric energy ${\mathcal E}_{pot}$ (semi-log scale) for the Landau damping test using the Hamiltonian Strang splitting with 32x64x64 grid points and $\Delta t=0.05$. Top figure: nonlinear Landau damping ($\alpha=0.5$); bottom figure: linear Landau damping ($\alpha=0.01$). The analytic decay rate $\gamma\approx 0.066$ is shown for comparison in case of the linear Landau damping.} \label{fig:ld-ee} 
\end{figure}

In Figure \ref{fig:ld-ee}, the time history of the electric energy is shown for a Strang Hamiltonian splitting, 
using $32$ points in the spatial direction, $64$ points in each velocity direction and a time step $\Delta t=0.05$. 
Note that compared to the linear case ({\it i.e.}~$\alpha=0.01$, for example) the electric energy does only decay exponentially for a brief period of time before presenting an oscillating behavior for large times. 
The results is in a good agreement with results available in the literature (see \cite{manfredi}). 

\begin{figure}
	\includegraphics[width=12cm]{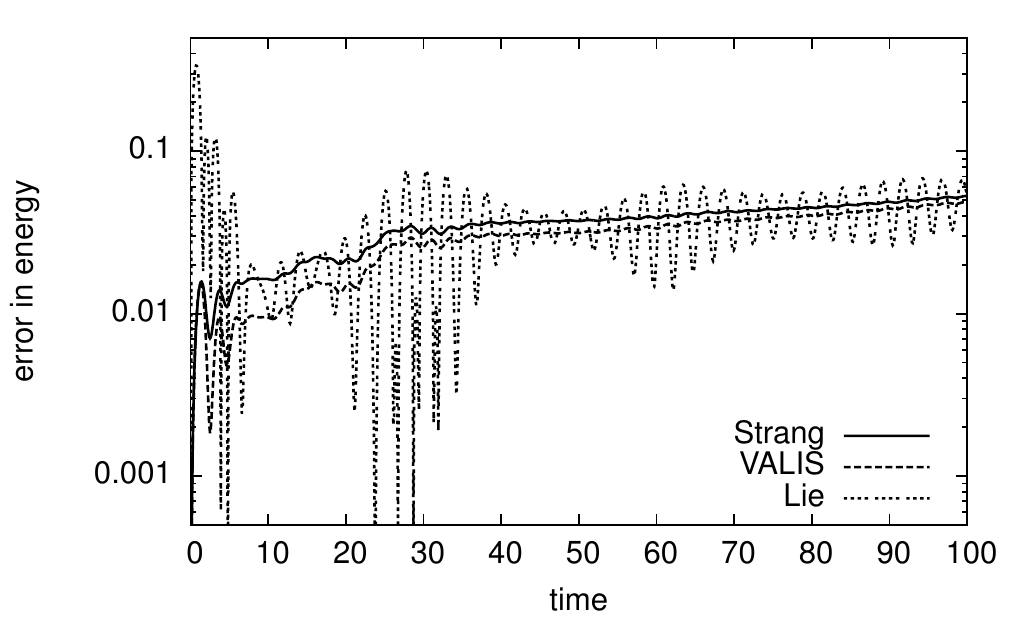}
 	\includegraphics[width=12cm]{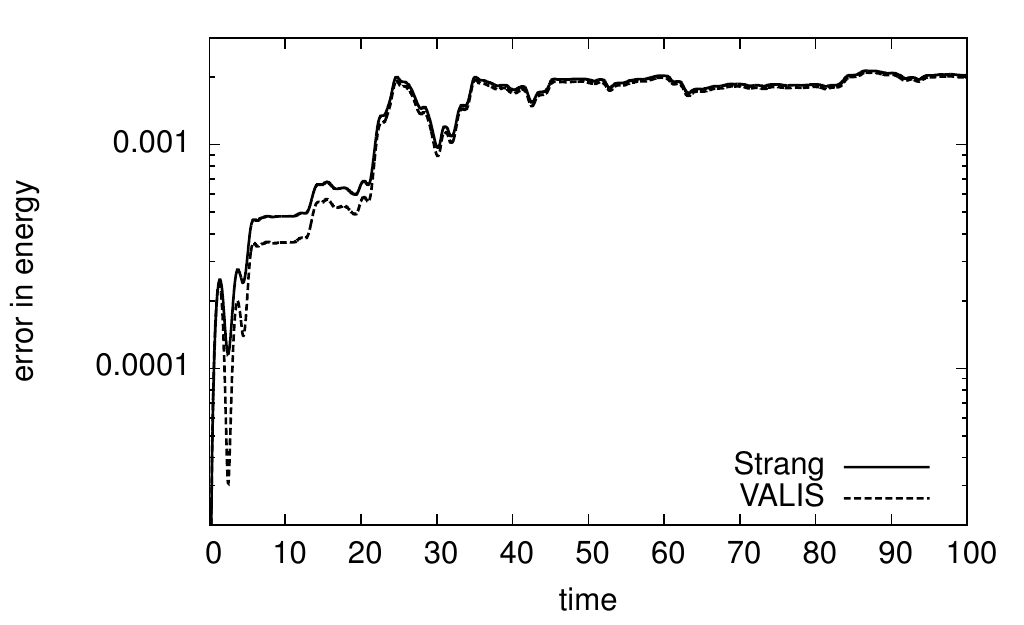}
	\caption{Time evolution of the error in energy (semi-log scale) for the Landau damping test $|{\mathcal E}_{tot}(t)-{\mathcal E}_{tot}(0)|$. 
	Top figure: 32x64x64 grid points and $\Delta t=0.1$; bottom figure: 256x256x256 grid points and $\Delta t=0.0125$. }
	\label{fig:landau-energy-error}
	\end{figure}

In Figure \ref{fig:landau-energy-error} the time history of the total energy error is plotted for a Lie and Strang Hamiltonian splitting (the scheme proposed in this paper) as well as the  VALIS scheme described in \cite{valis} (for more details see appendix \ref{sec:other-algorithms}). On the upper figure, 
we consider a phase space mesh of $32\times 64\times 64$ and the time step is $\Delta t=0.1$. 
We can first observe that the three methods behave well as the error in energy is less than $1\%$. 
However, the Lie splitting leads to oscillations which is not the case when a second order 
charge preserving method (Strang splitting or VALIS) is used. When the numerical parameters 
are refined (lower figure), these two methods  show virtually indistinguishable results  (the error is about 
$0.1 \%$) and we do not seem to experience a significant growth over the time period considered in this simulation.  
Let us remark that the charge conservation is very important with respect to energy conservation. If charge conservation is not exact, the accuracy of the algorithm is significantly reduced. This is more pronounced if we consider an increased number of grid points in space. The scheme proposed in \cite{mangeney} for that reason only performs satisfactory if a significantly smaller time step is chosen (as compared to Hamiltonian splittings and VALIS). This phenomenon can be observed in all the simulations conducted in this section and we will therefore omit it from the plots. We will, however, show a comparison for a reduced step size in the fully magnetized case (section \ref{sec:tsi}).

\subsection{Weibel instability}\label{sec:tt}
In this section, we will discuss an incarnation of the Weibel instability. This instability includes a genuine magnetic effect. Following \cite{pcp}, we impose the following initial particle density
\[
	f_0(x_1,v_1,v_2) = \frac{1}{\pi v_{\text{th}} \sqrt{T_r}} \e^{(v_{1}^2 + v_2^2/T_r)/v_{\text{th}}}(1+\alpha \cos(k x_1)), 
\]
where $x_1\in X=[0, 2\pi/k]$, $v\in \R^2$ and we have chosen \mbox{$v_{\text{th}}=0.02$}, $T_r=12$, $k=1.25$, and $\alpha=10^{-4}$. The electric field $E_1$ is initialized by solving Poisson's equation and $E_2$ is set to zero at $t=0$. The magnetic field is prescribe as 
	\[
		B_0(x_1) = \beta \cos(kx_1),
	\]
	where we have chosen $\beta=10^{-4}$.

\begin{figure}
		\includegraphics[width=12cm]{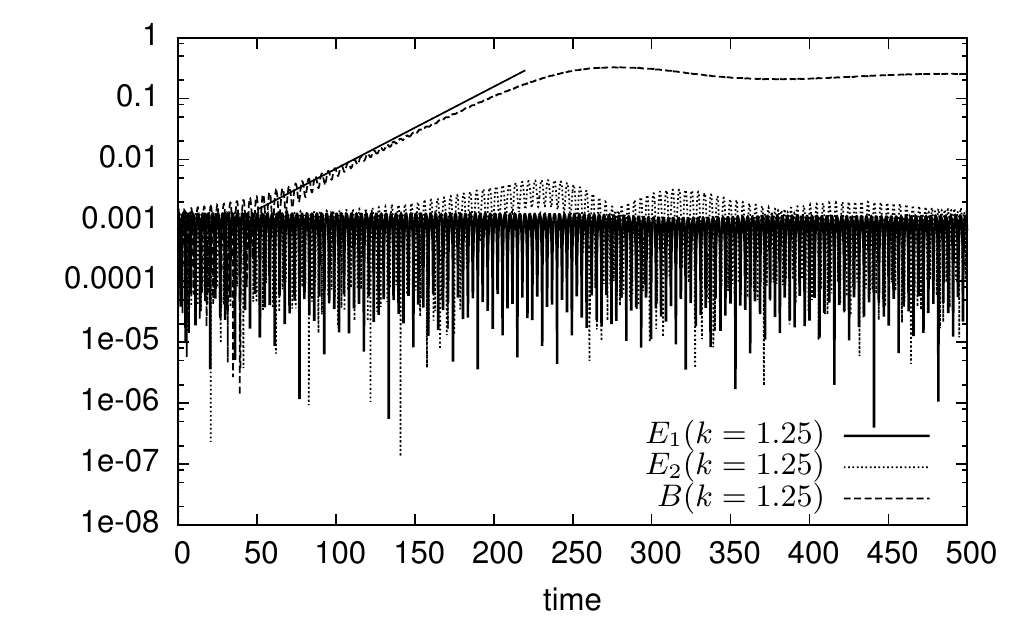}

		\includegraphics[width=12cm]{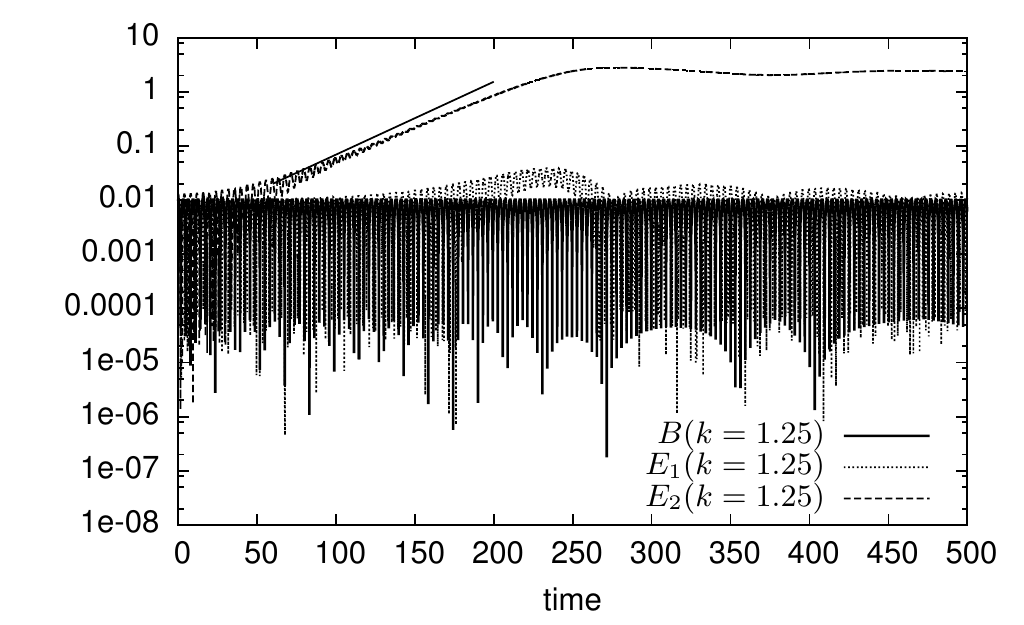}
		\caption{Time evolution of the $k=1.25$ Fourier mode of the electric and magnetic fields  (semi-log scale)
		for Weibel instability using the Hamiltonian Strang splitting method. 
		Top figure: 32x64x64 grid points and $\Delta t=0.2$; bottom figure: 256x256x256 
		grid points and $\Delta t=0.0125$.} \label{fig:tt-evolution}
	\end{figure}

	From Figure \ref{fig:tt-evolution} we see that the electric field in the $y$-direction 
does show an exponential growth in magnitude. The (analytically derived) growth rate is proportional to $k/\sqrt{1+k^2}$ and, for the parameters, chosen here is approximately equal to $0.031$. We observe that, similarly  to \cite{valis}, the agreement to the theory is as well as one would expect given that the analytical growth rate is derived using the linearized Vlasov--Maxwell equations. 
 Let us also note that even though the qualitatively features of the solution are well resolved on a grid with 32x64x64 points, quantitatively the result is wrong by at least an order of magnitude. 

	As in the previous example, we observe that the scheme proposed by Mangeney does not produce accurate results with respect to energy conservation, except if  the chosen step size is  significantly  smaller than what is necessary in case of the VALIS
	or the Hamiltonian splitting scheme (several orders of magnitude). However, compared to the simulation presented in the previous section, for sufficiently long integration times, we observe a better energy conservation for the Hamiltonian Strang splitting scheme as compared to the VALIS implementation (which is a scheme of second order as well).

	\begin{figure}
		\includegraphics[width=12cm]{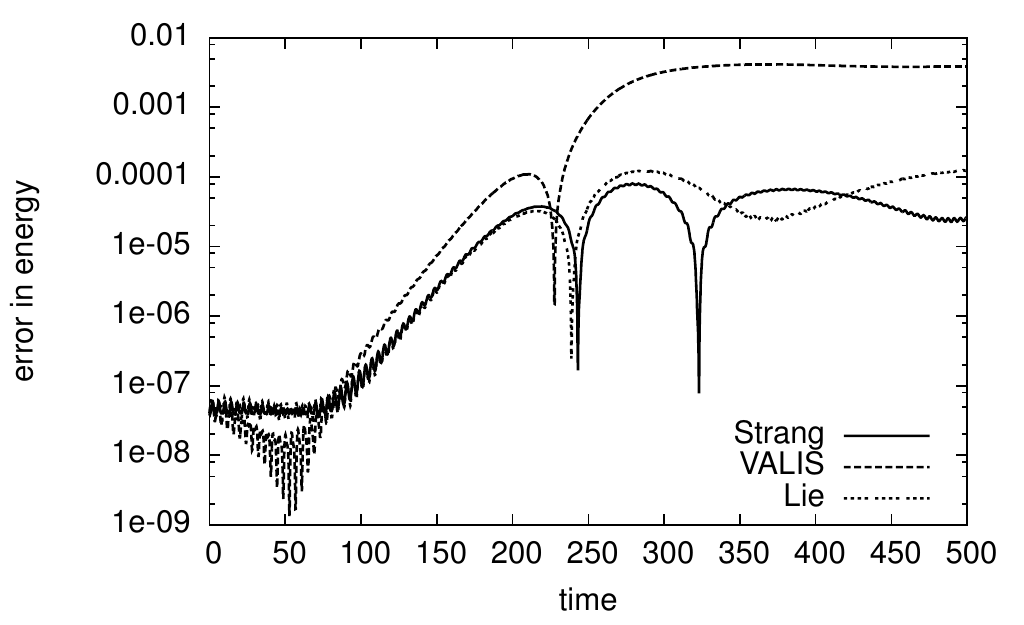}

		\includegraphics[width=12cm]{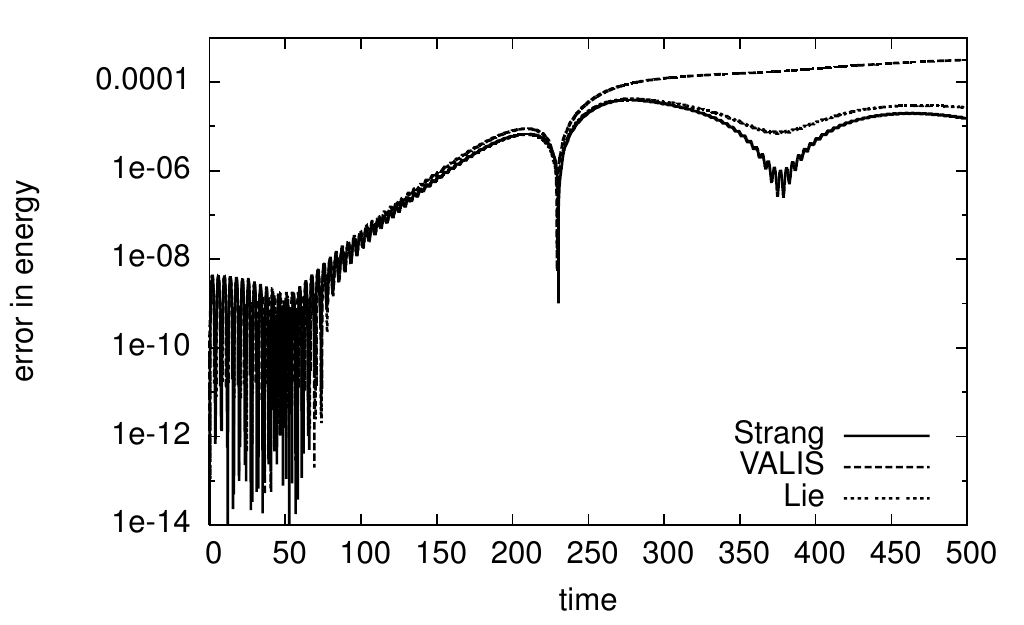}
		\caption{Weibel test: time evolution of the total energy error 
		$|{\mathcal E}(t)-{\mathcal E}(0)|$  (semi-log scale). 
		Top figure: 32x64x64 grid points and $\Delta t=0.2$; bottom figure: 256x256x256 grid points 
		and $\Delta t=0.0125$.} \label{fig:tt-energy}
	\end{figure}
	
\begin{figure}
		\includegraphics[width=12cm]{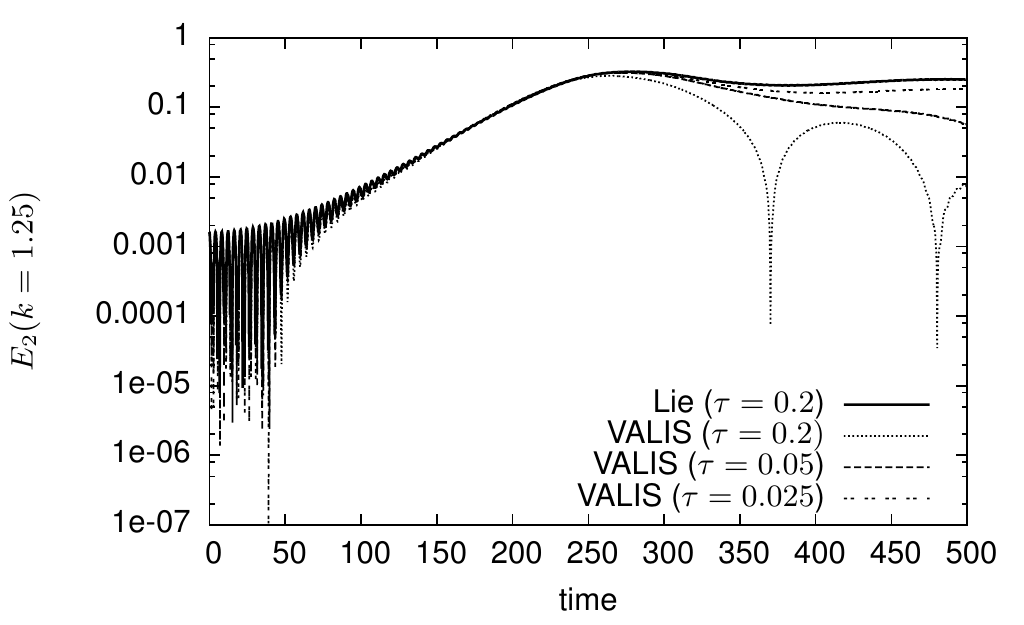}

		\caption{Weibel test: time evolution of the first mode of $E_2$ (semi-log scale). Comparison between VALIS, 
		and Strang splitting scheme using 32x64x64 grid points. 
		}\label{fig:tt-valis}
	\end{figure}

	Finally, let us investigate the behavior of the VALIS scheme more closely. 
	As can be seen from Figure \ref{fig:tt-valis}, if a time step close to the CFL 
	condition is chosen we can observe (unphysical)  oscillations in the solution. 
	Let us duly note that such oscillations are not present in our Hamiltonian 
	splitting approach. 
	In fact, we can verify that these oscillations disappear for the VALIS scheme when the time step size 
	is sufficiently small. This behavior implies that for the 
	problem under considered in this section, the VALIS scheme has a significantly 
	higher computational cost compared to the Hamiltonian Strang splitting 
	scheme proposed in this paper. 
	
\subsection{Two-stream instability} \label{sec:tsi}
As the final example, we will consider the two-stream instability. That is, the initial particle density is constructed from two Gaussian velocity profiles which are different with respect to their average velocity
	\[
		f_0(x_1,v_1,v_2) = \frac{1}{2 \pi \beta} \e^{-v_2^2/\beta}(\e^{-(v_1-0.2)^2/\beta} + \e^{-(v_1+0.2)^2/\beta}),
	\]
	where $x_1\in X=[0, 2\pi]$, $v\in \R^2$ and $\beta=2\cdot10^{-3}$. The resulting instability is usually driven by a perturbation in the particle density. However, as we are more interested in magnetic effects we will drive the instability by a perturbation in the magnetic field only
	\[
		B_0(x_1) = \alpha \sin(x_1),
	\]
	where we have chosen $\alpha=10^{-3}$.

\begin{figure}
		\includegraphics[width=12cm]{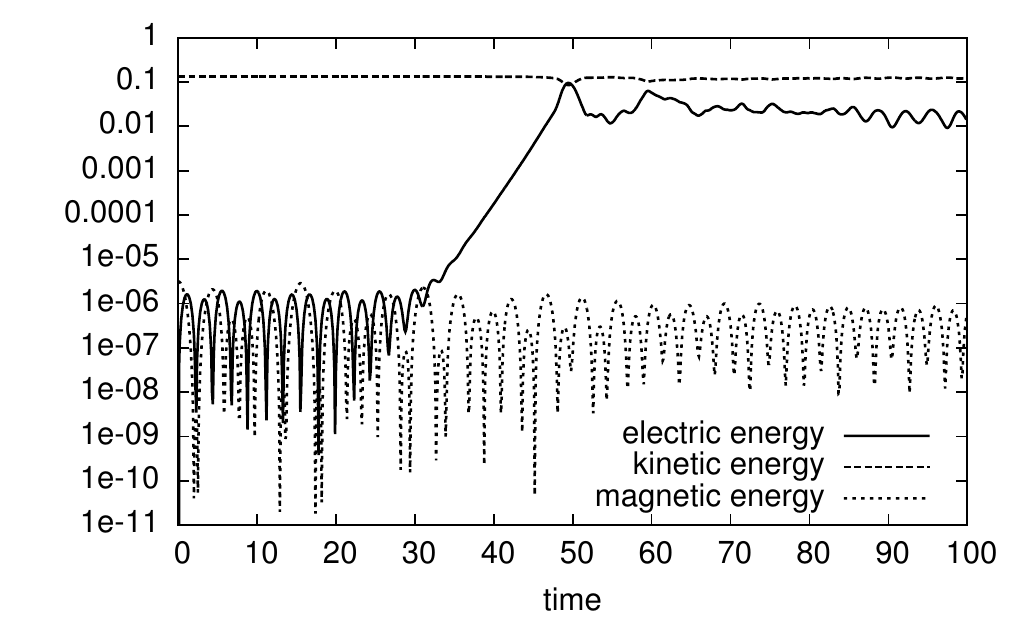}
		\caption{Time evolution of the electric, potential and total energies (semi-log scale) for the two-stream instability test using the Strang splitting scheme with 32x64x64 grid points and a time step of $\Delta t=0.1$.} \label{fig:tsi-evolution}
	\end{figure}

	The result of the simulation is shown in Figure \ref{fig:tsi-evolution}. Even though we drive the instability by a magnetic perturbation only, we observe, after an initial oscillatory regime, an exponential increase and a saturation in the electric field. 
	
\begin{figure}
		\includegraphics[width=12cm]{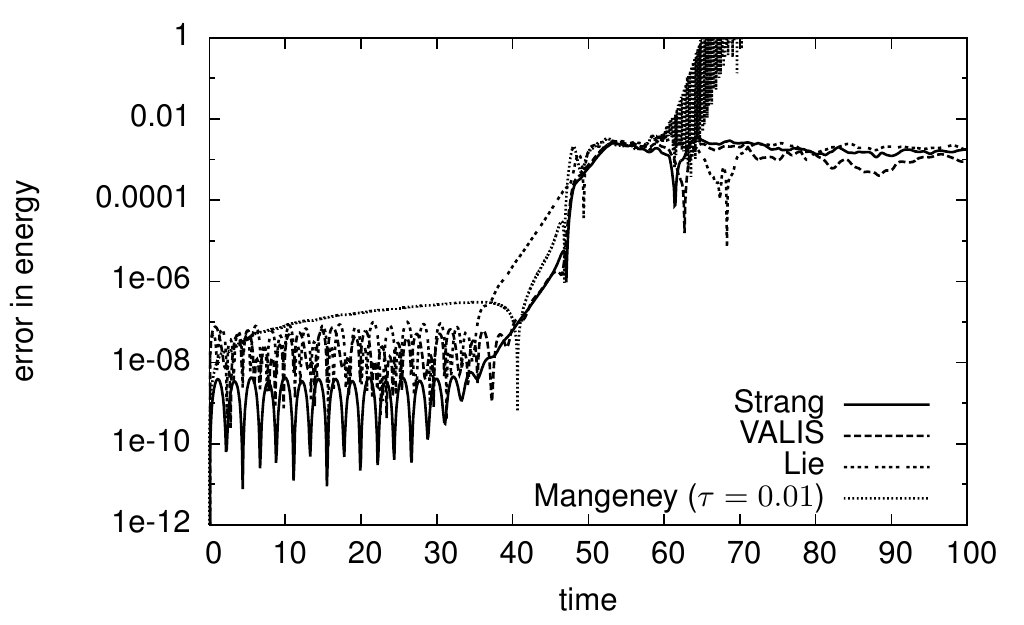}

		\includegraphics[width=12cm]{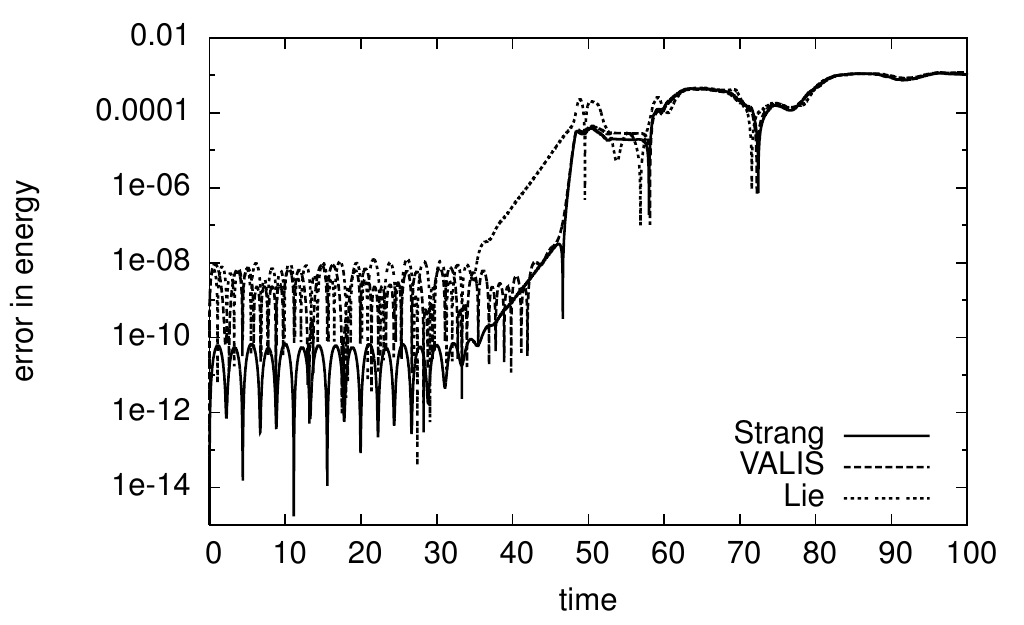}
		\caption{Time evolution of the error in energy $|{\mathcal E}_{tot}(t)-{\mathcal E}_{tot}(0)|$ (semi-log scale) of the two-stream instability test. 
		Top figure: 32x64x64 grid points and $\Delta t=0.1$ for VALIS, Strang and Lie, and 
		$\Delta t=0.01$ for Mangeney; 
		bottom figure: 256x256x256 grid points and $\Delta t=0.0125$.} \label{fig:tsi-energy}
	\end{figure}

	In Figure \ref{fig:tsi-energy} the error in energy for different schemes is shown. 
	Note that for the scheme proposed by Mangeney, even for a relatively coarse space discretization, 
	the step size has to be reduced by at least an order of magnitude in order to get comparable 
	results to the VALIS and Hamiltonian Strang splitting schemes. In addition, we observe that, 
	especially for a sufficiently fine space discretization, the error in energy for the VALIS and our 
	Hamiltonian Strang splitting scheme is almost indistinguishable.

\subsection{Order}

To conclude the discussion we consider the order achieved by Hamiltonian Lie and Strang splitting schemes as well as for the triple jump scheme. The order of the triple jump scheme is of particular interest as the results provided in this section demonstrate that we can in fact use composition to obtain a method that is of order four. The numerical results for the Weibel instability are shown in Table \ref{tab:order-64} for a space discretization using $64$ grid points and in Table \ref{tab:order-256} for a space discretization using $256$ grid points. In both cases we compute the $l^1$ error of the field variables ({\it i.e.} the two components of the electric field and the single component of the magnetic field) by using a reference solution for which a sufficiently small time step has been chosen.

\begin{table}
	\centering \textbf{Lie splitting scheme}

	\begin{tabular}{lll}
		stepsize & $l^1$ error & order \\
		\hline
		0.2 	&	1.881556551139132e-06& \\
		0.1 	&	8.769892305391039e-07&	 1.1013	\\
		0.05 	&	4.237597343971976e-07&	 1.04931	\\
		0.025 	&	2.083331679902264e-07& 	 1.02435	\\
		0.0125 	&	1.033058984634212e-07& 	 1.01197	\\
		0.00625 &	5.144678908873264e-08&	 1.00577	\\
	\end{tabular}

	\vspace{0.2cm}
	\textbf{Strang splitting scheme}

	\begin{tabular}{lll}
		stepsize & $l^1$ error & order \\
		\hline
		0.2		& 2.124324456920513e-07& \\
		0.1		& 5.27651798019974e-08 &	2.00935\\
		0.05	& 1.315400004930548e-08& 	2.00408\\
		0.025	& 3.276015657686228e-09& 	2.00549\\
		0.0125	& 8.068242472080718e-10& 	2.02162\\
		0.00625	& 1.86093775964047e-10 &	2.11622\\
	\end{tabular}

	\vspace{0.2cm}
	\textbf{Triple jump splitting scheme}

	\begin{tabular}{lll}
		stepsize & $l^1$ error & order \\
		\hline
		0.2 	&	1.321445314225445e-08& \\
		0.1 	&	8.095000243664432e-10& 	4.02894\\
		0.05 	&	5.162111972126581e-11&	3.971\\
		0.025 	&	2.031275496535852e-12&	4.6675\\
		0.0125 	&	4.297134131610342e-13& 	2.24094 \\
	\end{tabular}
	\vspace{0.2cm}
	\caption{Order determination for the Hamilton Lie, Strang, and triple jump splitting scheme in case of the Weibel instability introduced in section \ref{sec:tt}. In these simulations $64$ grid points are used in each dimension and the system is integrated up to final time $t=1$.} \label{tab:order-64}
\end{table}

\begin{table}
	\centering \textbf{Lie splitting scheme}

	\begin{tabular}{lll}
		stepsize & $l^1$ error & order \\
		\hline
		0.2  &	1.823383817066824e-06 &        \\ 
		0.1	 &	 8.465392041182712e-07& 1.10697\\
		0.05 &	 4.082184424389359e-07& 1.05224\\
		0.025&	 2.004799316967084e-07& 1.02588\\
	\end{tabular}

	\vspace{0.2cm}
	\textbf{Strang splitting scheme}

	\begin{tabular}{lll}
		stepsize & $l^1$ error & order \\
		\hline
		0.2		& 2.271250317272867e-07&             \\
		0.1		& 5.640193594261031e-08&	 2.00967 \\
		0.05	& 1.407771489982866e-08&	 2.00233 \\
		0.025	& 3.51820025143404e-09 &  2.0005     \\
	\end{tabular}

	\vspace{0.2cm}
	\textbf{Triple jump splitting scheme}

	\begin{tabular}{lll}
		stepsize & $l^1$ error & order \\
		\hline
		0.05 & 5.922461122797243e-11 & \\
		0.025 &3.371093950592525e-12 & 4.13491 \\
	\end{tabular}
	\vspace{0.2cm}
	\caption{Order determination for the Hamilton Lie, Strang, and triple jump splitting scheme in case of the Weibel instability introduced in section \ref{sec:tt}. In these simulations $256$ grid points are used in each dimension and the system is integrated up to final time $t=1$.} \label{tab:order-256}
\end{table}

\section{Appendix: Other algorithms} \label{sec:other-algorithms}
We present in this section some algorithms found in the literature \cite{valis}. The algorithms only differ from the 
computation of the Maxwell equations. The computation of the distribution function through the Vlasov equation 
is always performed using a Strang splitting 
\begin{itemize}
\item $\hat{f}^{\star} = \hat{f}^n \exp(-ik v_1 \Delta t/2)$ solution of $\partial_t f +v_1\partial_{x_1}f = 0$  on $\Delta t/2$, 
\item compute $f^{\star\star}$ by solving  $\partial_t f +(E+B {\mathcal J}v)\cdot \nabla_v f  = 0$  on $\Delta t$, 
using a splitting procedure,  
\item $\hat{f}^{n+1} = \hat{f}^{\star\star} \exp(-ik v_1 \Delta t/2)$ solution of  $\partial_t f +v_1\partial_{x_1}f = 0$  on $\Delta t/2$.  
\end{itemize}

\subsection{Predictor-Corrector}
\label{algo-mangeney}
At time $t^n$, the electric field $E^{n}$ and the magnetic field $B^{n-1/2}$ are known, together with $f^n$. 

\begin{itemize}
\item compute the current $J_{1, 2}^n(x_1)=\int_{\R^2} v_{1, 2} f^n(x_1, v)\dd v$
\item compute $B^{n+1/2}=B^n-\Delta t/2 \; \partial_{x_1}E_2^{n}$
\item compute $E_2^{n+1/2}=E_2^n-\Delta t/2 \; \partial_{x_1}(B^{n-1/2}+B^{n+1/2})/2 -\Delta t/2 J_2^n$ 
\item compute $E_1^{n+1/2}=E_1^n -\Delta t/2 \; J_1^n$ 
\item compute $f^{n+1}$ with Strang splitting, the fields $(E^{n+1/2}, B^{n+1/2})$ being known at time $t^{n+1/2}$, 
\item compute the current $J_{1, 2}^{n+1}=\int_{\R^2} v_{1, 2} f^{n+1}(x_1, v)\dd v$ 
\item compute $E_2^{n+1}=E_2^n-\Delta t \partial_{x_1}B^{n+1/2} -\Delta t/2 (J_2^n+J_2^{n+1})$ 
\item compute $E_1^{n+1}=E_1^n -\Delta t/2 (J_1^n+J_1^{n+1})$ 
\end{itemize}

\subsection{Valis}
\label{algo-valis}
At time $t^n$, the electric field $E^{n}$ and the magnetic field $B^{n-1/2}$ are known, together with $f^n$. 

Note that this algorithm conserves the charge.

\begin{itemize}
\item compute the current $J_{1, 2}^n(x_1)=\int_{\R^2} v_{1, 2} f^n(x_1, v)\dd v$
\item compute $B^{n+1/2}=B^n-\Delta t/2 \; \partial_{x_1}E_2^{n}$
\item compute $E_2^{n+1/2}=E_2^n-\Delta t/2 \; \partial_{x_1}(B^{n-1/2}+B^{n+1/2})/2 -\Delta t/2 J_2^n$ 
\item compute $E_1^{n+1/2}=E_1^n -\Delta t/2 \; J_1^n$ 
\item compute $f^{n+1}$ and the current 
\begin{itemize}
\item from $\hat{f}^{\star}$, compute and store \\% = \hat{f}^n \exp(-ik v_1 \Delta t/2)$ and \\ store 
$\hat{J}^\star_{1, 2}= \int_{\R^2} v_{1,2} \hat{f}^\star \dd v=\int_{\R^2} \hat{f}^n (\exp(-ik v_{1,2} \Delta t/2)-1)/(ik \Delta t/2)\dd v$
\item compute $f^{\star\star}$ by solving the (conservative) advection in $v_1, v_2$, 
the fields being known at time $t^{n+1/2}$, 
\item from $\hat{f}^{n+1}$, compute and store \\ %
$\hat{J}^{\star\star}_{1,2}= \int_{\R^2} v_{1,2} \hat{f}^{n+1} \dd v = \int_{\R^2} \hat{f}^{\star\star} (\exp(-ik v_{1,2} \Delta t/2)-1)/(ik \Delta t/2)\dd v$
\end{itemize}
\item compute the current $\hat{J}_{1,2}^{n+1/2} = (\hat{J}^{\star}_{1,2}+\hat{J}^{\star\star}_{1,2})/2$,  
\item compute $E_2^{n+1}=E_2^n-\Delta t \partial_{x_1}B^{n+1/2} -\Delta t J_2^{n+1/2}$ 
\item compute $E_1^{n+1}=E_1^n -\Delta t J_1^{n+1/2} $ 
\end{itemize}

In \cite{valis}, the current is extracted from the two $x_1$-advection steps. Using Fourier transforms, 
it  is written 
$$
\hat{f}^{\star} = \hat{f}^{n}  \exp(-ik v_1 \Delta t/2) %
  =  \hat{f}^{n} -ik\frac{\Delta t}{2} \hat{f}^{n}  \left[ \exp(-ik v_1 \Delta t/2) - 1 \right]/(ik \Delta t/2) , \nonumber\\
$$
which becomes, after integration in $v\in\R^2$, $\hat{\rho}^{\star} =\hat{\rho}^n -ik\frac{\Delta t}{2} \hat{J}_1^\star$ with 
$$
 \hat{J}_1^\star=  \int_{\R^2}  \hat{f}^{n}  \left[ \exp(-ik v_1 \Delta t/2) - 1 \right]/(ik \Delta t/2) \dd v.  
$$
Similarly the second part of the current  is denoted by 
$$
\hat{J}_1^{\star\star}=\int_{\R^2}  \hat{f}^{\star\star}  \left[ \exp(-ik v_1 \Delta t/2) - 1 \right]/(ik \Delta t/2) \dd v. 
$$  
Gathering these two steps, we can compute the density at time $t^{n+1}$ by integrating in $v$ each step of the splitting for $f$ 
\begin{eqnarray*}
\hat{\rho}^{n+1} &=& \hat{\rho}^{\star\star}  -ik\frac{\Delta t}{2} \hat{J}_1^{\star\star}\nonumber\\
                 &=& \hat{\rho}^{\star}  -ik\frac{\Delta t}{2} \hat{J}_1^{\star\star}\;\;\;  \mbox{(since the advection in $v$ is conservative)}\nonumber\\
&=& \hat{\rho}^{n} -ik\frac{\Delta t}{2} \hat{J}_1^{\star} -ik\frac{\Delta t}{2} \hat{J}_1^{\star\star}	\nonumber\\
&=& \hat{\rho}^{n} -ik\frac{\Delta t}{2} \left[ \hat{J}_1^{\star} +\hat{J}_1^{\star\star}\right]	\nonumber\\
&=& \hat{\rho}^{n} -ik \Delta t \hat{J}_1^{n+1/2}	\nonumber\\
\end{eqnarray*}
Assuming the Poisson equation satisfied at time $t^n$, we get $ik \hat{E}^n_1 = \hat{\rho}^n$ and since the current is used $\hat{J}_1^{n+1/2}$ 
to advance the Amp\`ere equation $\hat{E}_1^{n+1} = \hat{E}_1^{n} - \Delta t \hat{J}_1^{n+1/2}$, we have 
$$
\hat{\rho}^{n+1} = ik \hat{E}^n_1+ ik (\hat{E}_1^{n+1} - \hat{E}_1^{n}) =  ik \hat{E}_1^{n+1}, 
$$
which is exactly the Poisson equation at time $t^{n+1}$. 
Note that this can be generalized for finite volume method following the ideas of \cite{valis}. 

\section{Conclusion}
In this work, a new time splitting is proposed for the numerical solution of the  Vlasov--Maxwell system. 
This splitting is based on a decomposition of the Hamiltonian of the Vlasov--Maxwell system. 
Each step can be solved exactly in time so that the error is only due to the splitting procedure 
(first order for Lie, second order for Strang, fourth order for the triple jump scheme). In addition to the 
fact that methods of arbitrary high order can be constructed by composition, we prove 
that when one uses a spectral method for the spatial discretization, the charge is exactly preserved. 
Several numerical tests show the very good behavior of the splitting compared to methods from the 
literature, especially regarding the conservation of total energy.  

Several further perspectives of this work can be envisaged. First, a generalization to finite volumes methods for the spatial discretization have to be performed, taking care of the charge conservation property. We are confident this can be achieved by combining ideas of \cite{valis} and \cite{crous-respaud}. 
Second, a more detailed study of high order splittings (such as the triple Jump scheme) is necessary, since the use of small steps due to the CFL condition imposed by the Maxwell equations is a limiting factor. A substepping methodology would alleviate this problem. Finally, considering 
more realistic problems in $4$ or $5$ dimensions will require the parallelization of the code in the context of distributed memory architectures ({\it i.e.}~using MPI). Due to the fact that the VALIS method requires more temporary variables than the splitting method proposed in this paper, a more efficient MPI parallelization of our Hamiltonian splitting seems to be possible and this will constitute further research.

\end{document}